# Real-Time Optimal Power Flow under Wind Energy Penetration-Part II: Implementation

Erfan Mohagheghi, *Student Member*, *IEEE*, Aouss Gabash, *Member*, *IEEE*, Pu Li
Department of Simulation and Optimal Processes
Institute of Automation and Systems Engineering
Ilmenau University of Technology
Ilmenau, Germany
erfan.mohagheghi@tu-ilmenau.de, aouss.gabash@tu-ilmenau.de, pu.li@tu-ilmenau.de

*Abstract*—In this paper (Part II) we implement the prediction-updating approach developed in Part I to address fast changes in wind power generation when solving a complex real-time optimal power flow (RT-OPF) problem. The approach considers essential scenarios around forecasted wind power values in a moving prediction horizon (120 seconds). The individual optimal power flow problems corresponding to these scenarios are solved in parallel using a multi-processor server. Then the operation strategy is updated in a short sampling time (every 20 seconds) considering real wind power values. The RT-OPF problem is formulated considering both technical and economic aspects simultaneously. The RT-OPF is implemented on a 41-bus medium-voltage distribution network with two wind stations. The results show the benefits of the proposed approach and highlight further challenges of RT-OPF.

*Keywords—Energy prices; parallel computing; real-time optimal power flow (RT-OPF); wind power curtailment*

NOMENCLATURE

**Sets**
| | |
|---|---|
| $sd$ | Set of demand buses, i.e., {4,6,8,10,13,14, 22,23,25,27,30,31,34,36,37,41}. |
| $sw$ | Set of wind station (WS) buses, i.e., {2,16}. |
| $sl$ | Set of wind power levels, i.e., {1,2,3}. |

**Functions**
| | |
|---|---|
| $f$ | Total value of objective function. |
| $f_1$ | Total revenue from wind power. |
| $f_2$ | Total cost of active energy losses in the grid. |
| $f_3$ | Total cost of active energy at slack bus. |
| $f_4$ | Total cost of reactive energy at slack bus. |
| $f_P$ | Network active power function. |
| $f_Q$ | Network reactive power function. |

**Parameters**
| | |
|---|---|
| $N$ | Total number of buses. |
| $P_d(i,t_{120s})$ | Active power demand at bus $i \in sd$ in prediction horizon $t_{120s}$. |
| $Price_p(t_{120s})$ | Price for active energy in prediction horizon $t_{120s}$. |
| $Price_q(t_{120s})$ | Price for reactive energy in prediction horizon $t_{120s}$. |
| $P_w(i,t_{120s})$ | Active power of WS at bus $i \in sw$ in prediction horizon $t_{120s}$. |
| $P_{W.r}(i)$ | Rated installed wind power at bus $i \in sw$. |
| $P_{w.A}(i,t_{20s})$ | Actual wind power of WS at bus $i \in sw$ in update interval $t_{20s}$. |
| $P_{w.H.\sigma}(i,t_{120s})$ | Wind power higher than forecasted of WS at bus $i \in sw$ in prediction horizon $t_{120s}$ for level $\sigma \in sl$. |
| $P_{w.L.\sigma}(i,t_{120s})$ | Wind power lower than forecasted of WS at bus $i \in sw$ in prediction horizon $t_{120s}$ for level $\sigma \in sl$. |
| $P_{w.M}(i,t_{120s})$ | Mean (forecasted) wind power of WS at bus $i \in sw$ in prediction horizon $t_{120s}$. |
| $Q_d(i,t_{120s})$ | Reactive power demand at bus $i \in sd$ in prediction horizon $t_{120s}$. |
| $S_{l.max}(i,j)$ | Upper limit of apparent power flow of line between bus $i$ and $j$. |
| $S_{S.max}$ | Upper limit of apparent power at slack bus. |
| $t_{120s}$ | Prediction horizon, i.e., 120 seconds. |
| $t_{112s}$ | Reserved time for computing OPF problems, i.e., 112 seconds. |
| $t_{20s}$ | Update interval, i.e., 20 seconds. |
| $V_{max}(i)$ | Upper limit of voltage at bus $i$. |
| $V_{min}(i)$ | Lower limit of voltage at bus $i$. |
| $\Delta P_\sigma(i)$ | Wind power deviation at bus $i$ for level $\sigma \in sl$. |
| $\mu_d$ | Mean value for demand power. |
| $\mu_w$ | Mean value for wind power. |
| $\sigma_d$ | Standard deviation for demand power. |
| $\sigma_w$ | Standard deviation for wind power. |

**Variables**
| | |
|---|---|
| $P_{loss}(t_{120s})$ | Active power losses in the grid in prediction horizon $t_{120s}$. |
| $P_S(t_{120s})$ | Active power injected at slack bus in prediction horizon $t_{120s}$. |
| $Q_S(t_{120s})$ | Reactive power injected at slack bus in prediction horizon $t_{120s}$. |
| $S(i,j,t_{120s})$ | Apparent power flow from bus $i$ to $j$ in prediction horizon $t_{120s}$. |
| $V(i,t_{120s})$ | Voltage at bus $i$ ($i \neq 1$) in prediction horizon $t_{120s}$. |
| $\beta_{c.w}(i,t_{120s})$ | Curtailment factor of wind power for WS at bus $i \in sw$ in prediction horizon $t_{120s}$. |
| $\beta_{c.w}(i,t_{20s})$ | Curtailment factor of wind power for WS at bus $i \in sw$ in update interval $t_{20s}$. |

This work is supported by the Carl-Zeiss-Stiftung.

The final version of this paper has been published in the proceeding of 2016 IEEE 16th International Conference on Environment and Electrical Engineering (EEEIC). ©2016 IEEE. DOI: 10.1109/EEEIC.2016.7555465

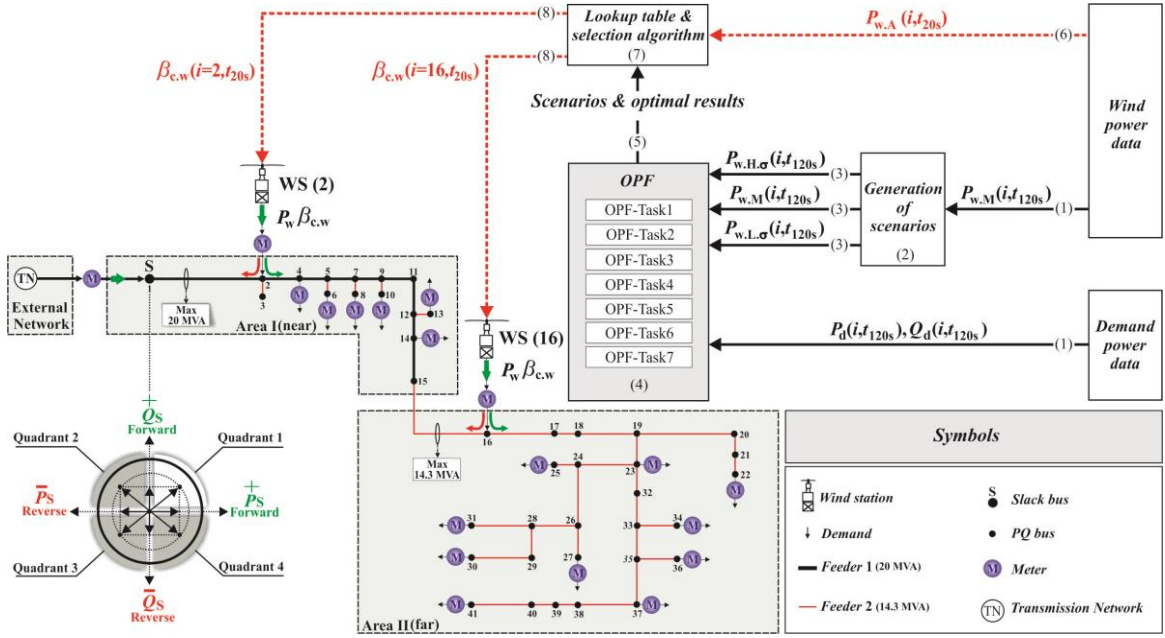

Fig. 1. Distribution network and the proposed RT-OPF framework.

## I. INTRODUCTION

The method of the real-time optimal power flow (RT-OPF) under wind energy penetration is presented in Part I of our paper [1]. For its implementation, we consider a distribution network (DN) shown in Fig. 1. It is assumed the two wind stations (WSs) have the same installed capacities at different locations in a DN, i.e., they are subject to different wind speeds. It is aimed in this paper to implement the RT-OPF for the network to achieve the following benefits:

- Considering both technical and economic issues in the formulation of the OPF problems;
- Ensuring not only optimal but also feasible operations for the network;
- Ensuring low computation time to solve the high-dimensional and highly complicated NLP-OPF problem in *real-time*.

## II. RT-OPF PROBLEM FORMULATION

The major task of the proposed RT-OPF is to solve the OPF problems corresponding to the essential scenarios around the forecasted wind power. The solution provides a lookup-table for the prediction horizon. In this work, we formulate each OPF problem for a prediction horizon with 120 seconds as follows:

$$\max_{\beta_{c.w}} \ f = f_1 - f_2 - f_3 - f_4 \quad (1)$$

$$f_1 = Price_p(t_{120s}) \sum_{\substack{i=1 \\ i \in sw}}^{N} P_w(i,t_{120s}) \beta_{c.w}(i,t_{120s}) \quad (2)$$

$$f_2 = Price_p(t_{120s}) \, P_{loss}(t_{120s}) \quad (3)$$

$$f_3 = Price_p(t_{120s}) \, P_S(t_{120s}) \quad (4)$$

$$f_4 = Price_q(t_{120s}) \, Q_S(t_{120s}) \quad (5)$$

where $t_{120s} = 120s$ is the prediction horizon. It is aimed in (1) to maximize the total revenue from the wind power $f_1$, and meanwhile to minimize the total costs of the active energy losses in the grid $f_2$, the cost of the active energy at slack bus $f_3$, and the cost of the reactive energy at slack bus $f_4$, respectively. Here, $Price_p(t_{120s})$ is the active energy price, $Price_q(t_{120s})$ the reactive energy price, and $P_{loss}(t_{120s})$ the active power losses in the grid in the time horizon $t_{120s}$. $N$ is the total number of buses, $P_w(i,t_{120s})$ is the active power of WS at bus $i$, while $sw$ stands for the set of WSs. $P_S(t_{120s})$ and $Q_S(t_{120s})$ are the active and reactive power injected at slack bus in the time horizon $t_{120s}$, respectively. The control variable of each WS is the curtailment factor of wind power for WS at bus $i$ i.e., $\beta_{c.w}(i,t_{120s})$.

The objective function expressed as Eqs. (1)-(5) is subject to the following equations as equality constraints:

$$f_P + P_d(i,t_{120s}) - P_w(i,t_{120s})\beta_{c.w}(i,t_{120s}) - P_S(t_{120s}) = 0, \ i \in N \quad (6)$$

$$f_Q + Q_d(i,t_{120s}) - Q_S(t_{120s}) = 0, \ i \in N \quad (7)$$

where (6) and (7) are the network active and reactive power equations [2], in which $f_P$ and $f_Q$ are the network active and reactive power functions, respectively. $P_d(i,t_{120s})$ and $Q_d(i,t_{120s})$ are, respectively, active and reactive power demand at bus $i \in sd$ in the prediction horizon $t_{120s}$. In addition, the following inequality constraints should be held which are the active and reactive bounds at slack bus

$$\left(P_S(t_{120s})\right)^2 + \left(Q_S(t_{120s})\right)^2 \leq \left(S_{S.max}\right)^2 \quad (8)$$

$$0 \leq P_S(t_{120s}) \leq S_{S.max} \quad (9)$$

$$0 \leq Q_S(t_{120s}) \leq S_{S.\max} \quad (10)$$

voltage bounds of the PQ-buses

$$V_{\min}(i) \leq V(i,t_{120s}) \leq V_{\max}(i), \qquad i \in N \ (i \neq 1) \quad (11)$$

feeder sections limits

$$S(i,j,t_{120s}) \leq S_{l.\max}(i,j), \qquad i \in N(i \neq j) \quad (12)$$

and the limits of the curtailment factors

$$0 \leq \beta_{c.w}(i,t_{120s}) \leq 1. \quad (13)$$

Here, $V(i,t_{120s})$ is the voltage at bus $i$ in the prediction horizon $t_{120s}$, $V_{\max}(i)$ and $V_{\min}(i)$ are respectively the upper and lower limits. $S_{S.\max}$ is the upper limit of apparent power at slack bus, and $S_{l.\max}(i,j)$ is the upper limit of apparent power flow of line between bus $i$ and $j$. The whole equality and inequality equations are coded in the GAMS framework and the resulting problem is solved using the NLP solver CONOPT3 in GAMS [3]. The scenarios around forecasted wind power $P_{w.M}(i,t_{120s})$, i.e., the higher-side $P_{w.H.\sigma}(i,t_{120s})$ and lower-side $P_{w.L.\sigma}(i,t_{120s})$ levels for the prediction horizon are generated as follows

$$P_{w.H.\sigma}(i,t_{120s}) = P_{w.M}(i,t_{120s}) + \Delta P_\sigma(i), \qquad \sigma \in sl \quad (14)$$

$$P_{w.L.\sigma}(i,t_{120s}) = P_{w.M}(i,t_{120s}) - \Delta P_\sigma(i), \qquad \sigma \in sl \quad (15)$$

where $\sigma$ is the wind power level resulting in $\Delta P_3(i) = 1.5\Delta P_2(i) = 3\Delta P_1(i)$. Fig. 2 shows the generated levels leading to 49 scenarios (i.e., combinations of levels) for each prediction horizon. Since each of the scenarios is independent, parallel computing is implemented so that the computation time can be significantly reduced. The parallel computing is carried out on a server with 2 (physical processors) Intel Xeon X5690, CPU 3.47 GHz (6 cores, 12 threads) and 64 GB RAM. As shown in Fig. 1 and described in Part I [1], one of the scenarios will be selected based on the actual wind power $P_{w.A}(i,t_{20s})$ in the current time interval (i.e., 20 seconds) and the corresponding curtailment factors ($\beta_{c.w}(i,t_{20s})$) are applied to the grid.

## III. CASE STUDIES

The network under consideration is taken from [4]-[6]. It is a 41-bus, 27.6 kV typical rural distribution network as shown in Fig. 1, where two different feeder capacities are considered. The peak power demand is 16.25 MVA [2] and the substation rating is 20 MVA. Two WSs each with ($P_{W.r}$ = 10 MW) and unity power factors (PFs) are located at buses 2 and 16, as shown in Fig. 1. The active and reactive energy prices are adopted from [4], [5] and fixed with 1.67 \$/MW.$t_{120s}$ and 0.4 \$/Mvar.$t_{120s}$, respectively, and the bus number 1 is considered as the slack bus. To implement the proposed RT-OPF strategy, the forecasted demand and wind power profiles as well as the actually realized wind power profiles should be available.

The active $P_d$ and reactive $Q_d$ power demand are assumed to follow the hourly IEEE-RTS fall season's days [4]. However, inside each hour, the demand profiles are generated for every bus at every 120 seconds by adding normally distributed random numbers (with $\mu_d = 0$ and $\sigma_d = 0.01$) to their hourly values. The resulting demand trajectories for one day (24 hours) are shown in Fig. 3(a).

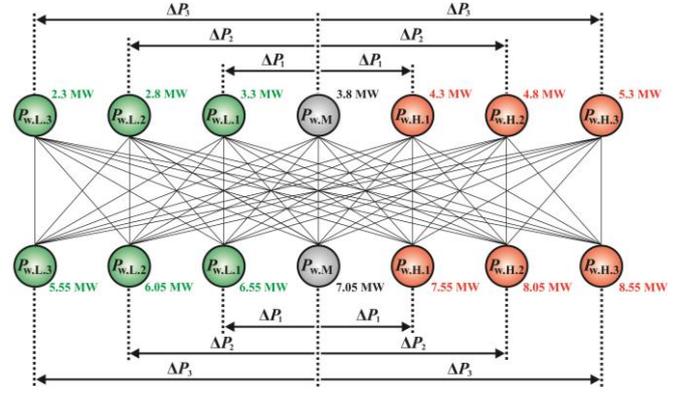

Fig. 2. Illustration of 49 scenarios for the two WSs (the values of the *first time horizon* are given).

In the practice, the forecasted and actual wind power profiles can be acquired from environmental data centers and online measurements, respectively. Here, we consider these profiles with random variations and generate them by simulation. The forecasted wind power profile for one day is taken from [4] and used to generate two different profiles for the two WSs. This is done by adding uniformly distributed random numbers ±15% for each WS at each hour to simulate the random behavior of wind energy. In addition, the forecasted wind power profiles are generated for every 120 seconds by adding normally distributed random numbers (with $\mu_w = 0$ and $\sigma_w = 0.1$) to their hourly values as shown in red-dashed curves in Fig. 3(b)-(c). Furthermore, the AWP for the two WSs are generated by adding uniformly distributed random numbers ±15% to the forecasted values for each WS at each 20 seconds as shown in blue-solid curves in Fig. 3(b)-(c). After these preparations, we run the RT-OPF described in [1] and Section 3 for one day as a demonstration. All trajectories are shown in Fig. 3 for 24 hours and for a selected time slot (8 minutes) in Fig. 4. Table I shows the lookup table for the first prediction horizon.

It can be seen in Table I that, for $t_{120s}1$, the active power at the slack bus is zero for all scenarios. This is because of the low active power demand of 6.65 MW and high active wind power generation. In addition, the high wind power generation leads to strong curtailment (i.e., low values of curtailment factors) to ensure a feasible operation. In contrast, the reactive power at the slack bus is about 2.3 Mvar for all scenarios in this time horizon, i.e., it will be imported from the upstream network. This is because of using unity PFs of the 2 WSs and the reactive power compensation of feeder capacitive susceptance [4] (the total reactive power demand is 2.46 Mvar in this time horizon). Therefore, the deviation between the forecasted and actual reactive power at the slack bus is insignificant, which can also be seen from Fig. 3(g) and 4(g).

To analyze the system performance in detail, 4 consecutive time horizons (i.e., 8 minutes) from Fig. 3 are selected and clearly presented in Fig. 4. Fig. 4(b) and 4(c) show the forecasted and actual profiles of the wind power based on which the optimal curtailment strategies are computed online for the two WSs, as shown in Fig. 4(d) and 4(e), respectively. As mentioned earlier, the optimization results from the proposed approach ensure feasible operations. As a result, the optimized strategy has some degree of conservatism due to the number of levels used for generating the scenarios. This can be seen from the resulting active power imported from the upstream network shown in Fig. 4(f), i.e., in most time the imported value is higher than the expected value. However, when the total wind power is much higher than the forecasted value but lower than the total demand plus the losses, the imported active power will be lower than the expected value, as seen from time step 1282 to 1285. It leads to higher total revenue in this time period, as shown in Fig. 4(h). Fig. 4(i) gives the computation time taken by the 7 processors for solving the optimization problems for 49 scenarios in real-time, showing that the computation time for each processor is less than the reserved time horizon $t_{112s}$ [1].

TABLE I. LOOKUP TABLE ARRANGEMENT FOR THE FIRST TIME HORIZON.

| Scenarios | | | | Optimal results of the scenarios for the prediction horizon | | | |
|---|---|---|---|---|---|---|---|
| No. | Scenario | $P_w(2,t_{120s}1)$ (MW) | $P_w(16,t_{120s}1)$ (MW) | $\beta_{c.w}(2,t_{120s}1)$ --- | $\beta_{c.w}(16,t_{120s}1)$ --- | $P_S(t_{120s}1)$ (MW) | $Q_S(t_{120s}1)$ (Mvar) |
| 1 | $P_{w.H.3}$-$P_{w.H.3}$ | 5.3 | 8.55 | 0.715 | 0.337 | 0 | 2.375 |
| 2 | $P_{w.H.3}$-$P_{w.H.2}$ | 5.3 | 8.05 | 0.715 | 0.358 | 0 | 2.375 |
| 3 | $P_{w.H.3}$-$P_{w.H.1}$ | 5.3 | 7.55 | 0.715 | 0.381 | 0 | 2.375 |
| 4 | $P_{w.H.3}$-$P_{w.M}$ | 5.3 | 7.05 | 0.715 | 0.408 | 0 | 2.375 |
| 5 | $P_{w.H.3}$-$P_{w.L.1}$ | 5.3 | 6.55 | 0.715 | 0.439 | 0 | 2.375 |
| 6 | $P_{w.H.3}$-$P_{w.L.2}$ | 5.3 | 6.05 | 0.715 | 0.476 | 0 | 2.375 |
| 7 | $P_{w.H.3}$-$P_{w.L.3}$ | 5.3 | 5.55 | 0.715 | 0.518 | 0 | 2.375 |
| ⋮ | ⋮ | ⋮ | ⋮ | ⋮ | ⋮ | ⋮ | ⋮ |
| 22 | $P_{w.M}$-$P_{w.H.3}$ | 3.8 | 8.55 | 0.997 | 0.337 | 0 | 2.375 |
| 23 | $P_{w.M}$-$P_{w.H.2}$ | 3.8 | 8.05 | 0.997 | 0.358 | 0 | 2.375 |
| 24 | $P_{w.M}$-$P_{w.H.1}$ | 3.8 | 7.55 | 0.997 | 0.381 | 0 | 2.375 |
| 25 | $P_{w.M}$-$P_{w.M}$ | 3.8 | 7.05 | 0.997 | 0.408 | 0 | 2.375 |
| 26 | $P_{w.M}$-$P_{w.L.1}$ | 3.8 | 6.55 | 0.997 | 0.439 | 0 | 2.375 |
| 27 | $P_{w.M}$-$P_{w.L.2}$ | 3.8 | 6.05 | 0.997 | 0.476 | 0 | 2.375 |
| 28 | $P_{w.M}$-$P_{w.L.3}$ | 3.8 | 5.55 | 0.997 | 0.518 | 0 | 2.375 |
| ⋮ | ⋮ | ⋮ | ⋮ | ⋮ | ⋮ | ⋮ | ⋮ |
| 43 | $P_{w.L.3}$-$P_{w.H.3}$ | 2.3 | 8.55 | 1 | 0.511 | 0 | 2.383 |
| 44 | $P_{w.L.3}$-$P_{w.H.2}$ | 2.3 | 8.05 | 1 | 0.543 | 0 | 2.383 |
| 45 | $P_{w.L.3}$-$P_{w.H.1}$ | 2.3 | 7.55 | 1 | 0.579 | 0 | 2.383 |
| 46 | $P_{w.L.3}$-$P_{w.M}$ | 2.3 | 7.05 | 1 | 0.62 | 0 | 2.383 |
| 47 | $P_{w.L.3}$-$P_{w.L.1}$ | 2.3 | 6.55 | 1 | 0.667 | 0 | 2.383 |
| 48 | $P_{w.L.3}$-$P_{w.L.2}$ | 2.3 | 6.05 | 1 | 0.722 | 0 | 2.383 |
| 49 | $P_{w.L.3}$-$P_{w.L.3}$ | 2.3 | 5.55 | 1 | 0.787 | 0 | 2.383 |

The implementation shows the applicability of the proposed RT-OPF, but the following open points are worthy of remark:

- In reality, reverse active and reactive power flows are likely to happen. However, in this work, reverse power flow is not allowed, in order to clearly show the impact of wind energy curtailments.
- It is assumed that for each prediction horizon (i.e., 120 seconds), the actual values of wind powers do not go beyond the scope of the levels in Fig. 2. Nevertheless, in reality, the actual wind power could be out of the levels which means more levels may be required.
- The width of all higher- and lower-side levels, $\Delta P_\sigma$, is considered to be the same. However, it could be more useful to determine the width of each level with respect to a probability density function.
- Although all solutions of the optimization problems in the case studies have been converged, for other problems, OPF could fail to reach the computation convergence. Hence a strategy to deal with convergence condition would be needed.
- Voltage at the slack bus is fixed with the amplitude of 1.0 pu and zero angle in this study. Since the assigned amplitude value may not be optimal for some scenarios, it would be more beneficial to consider it as a control variable to be optimized. This could make the problem even more complicated as large-scale *mixed-integer* NLP problems should be solved.

IV. CONCLUSIONS

In this paper, a RT-OPF was proposed and the method was implemented under changing wind energy penetration. A prediction-updating approach was used to deal with the random behavior of wind power generation. Parallel computing was implemented to speed up the computation. The implementation on a MV network demonstrates optimality and feasibility of the operation strategies of the RT-OPF framework. It was shown that the proposed RT-OPF ensures a short computation time which makes real-time optimization possible.

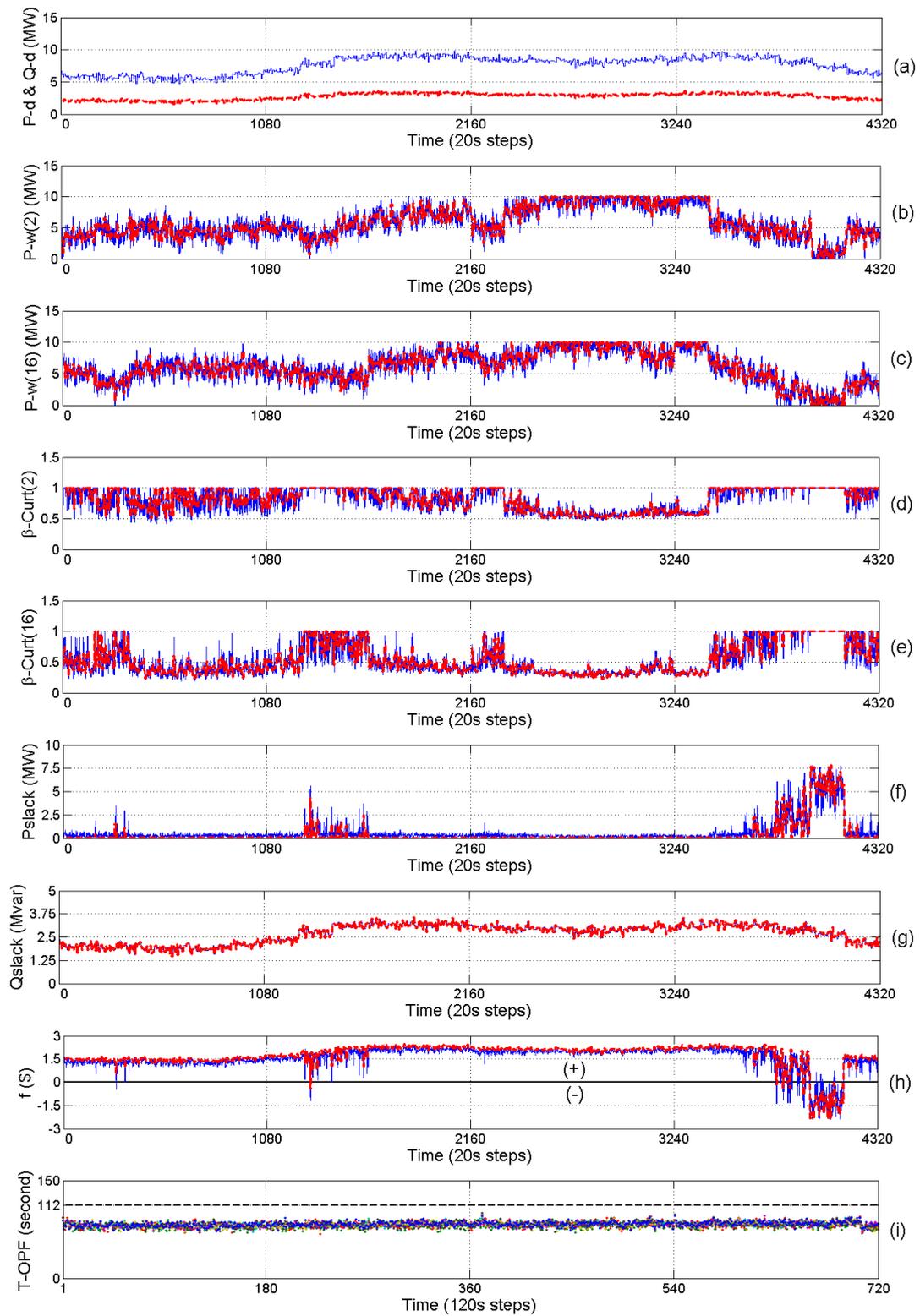

Fig. 3. Trajectories of one day. (a) Total active (blue-solid) and reactive (red-dashed) power demand. (b) Forecasted (red-dashed) and actual (blue-solid) wind power at bus 2. (c) Forecasted (red-dashed) and actual (blue-solid) wind power at bus 16. (d) Forecasted (red-dashed) and actual (blue-solid) curtailment factors at bus 2. (e) Forecasted (red-dashed) and actual (blue-solid) curtailment factors at bus 16. (f) Forecasted (red-dashed) and actual (blue-solid) slack bus active power. (g) Forecasted (red-dashed) and actual (blue-solid) slack bus reactive power. (h) Forecasted (red-dashed) and actual (blue-solid) total objective function value. (i) Computational time of the seven processors.

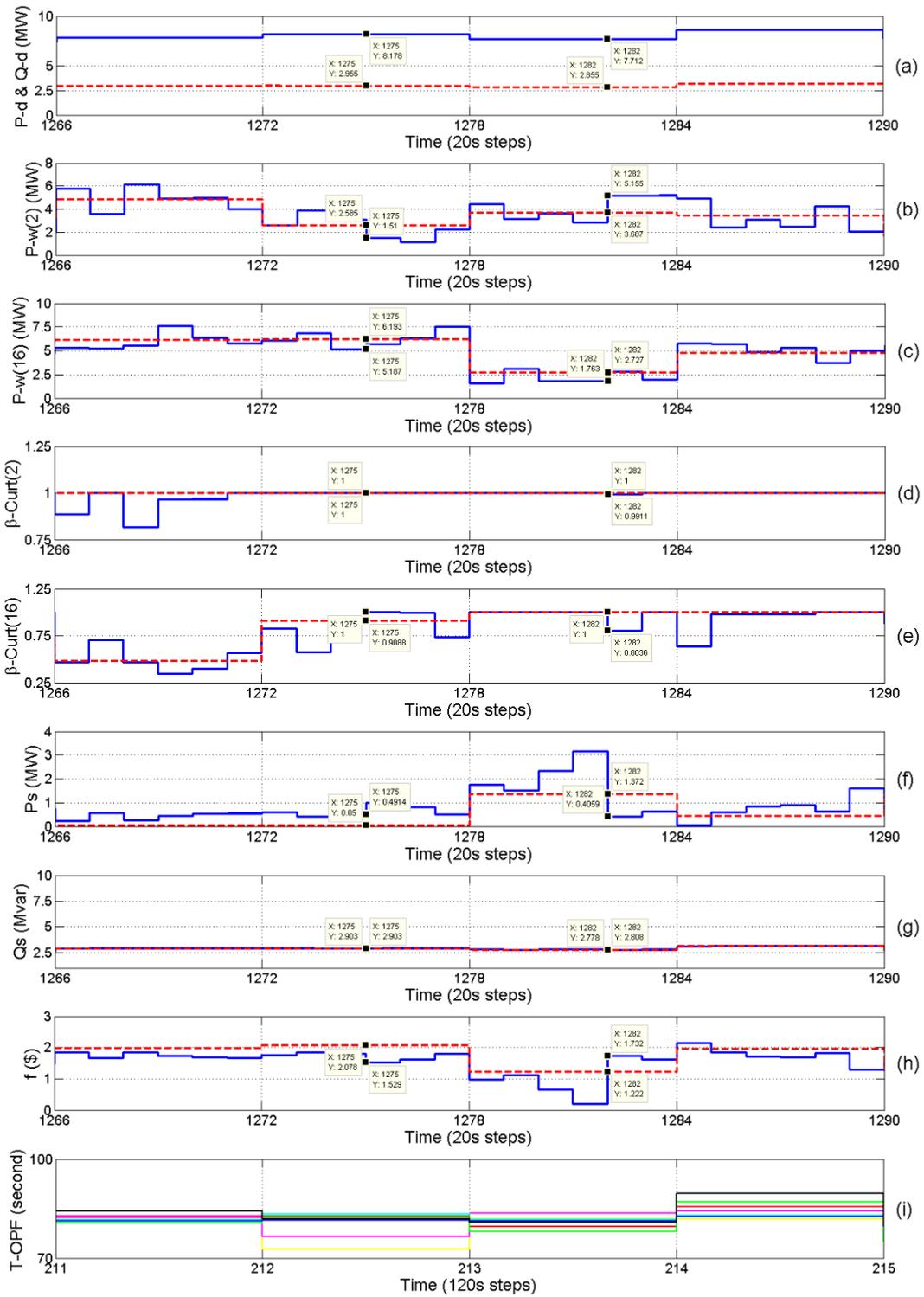

Fig. 4. The trajectories of the selected time slot (8 minutes) from Fig. 3.